\documentclass{article}
\usepackage{graphicx}
\usepackage{amsmath,amssymb}
\newcommand{\E}{\mathbb E}
\newcommand{\N}{\mathbb N}
\newcommand{\R}{\mathbb R}
\newcommand{\p}{\mathbb P}
\newcommand{\x}{\mathbf x}

\newcommand{\dd}{\text{\rm d}}
\newtheorem{corollary}{Corollary}
\newtheorem{proposition}{Proposition}
\newtheorem{definition}{Definition}
\begin{document}

\title{Planar segment processes with reference mark
distributions, modeling and estimation}
\author{Viktor Bene\v{s}, Jakub Ve\v{c}e\v{r}a, Milan Pultar\\\\
Charles University, Faculty of Mathematics and Physics \\ Department of Probability and Mathematical Statistics \\ Sokolovska 83, 18675 Praha 8, Czech Republic  }

\maketitle

\begin{abstract}The paper deals with planar segment processes given by a density with respect to the Poisson process. Parametric models involve reference distributions of directions and/or lengths of segments. These distributions generally do not coincide with the corresponding observed distributions. Statistical methods are presented which first estimate scalar parameters by known approaches and then the reference distribution is estimated non-parametrically. Besides a general theory we offer two models, first a Gibbs type segment process with reference directional distribution and secondly an inhomogeneous process with reference length distribution. The estimation is demonstrated in simulation studies where the variability of estimators is presented graphically.

\medskip

\noindent Keywords: Conditional intensity, segment process, semiparametric estimation

\medskip

\noindent AMS subject classification: 60D05, 60G55
\end{abstract}

\section{Introduction}\label{Introd}
The present research addresses an important problem in the statistics of spatial marked point processes given by a density with respect to the Poisson process. Observing a realization of spatial data which shall be fitted to such a model we first estimate the parameters by a method of point estimation. However, among the quantities to be estimated there may appear also the reference mark distribution which need not coincide with the observed mark distribution of the process. Both the scalar parameters and the reference mark distribution are needed e.g. when we try to simulate the model. This distribution can be also parametrized by a subjective choice of model, cf. \cite{RefBe}. In the present paper the main aim is to estimate the reference mark distribution non-parametrically, i.e. in total to use a semiparametric approach instead of a fully parametric one.

An early paper \cite{RefBa} mentions parameter estimation of a marked point process by means of the maximum pseudolikelihood method but the authors do not identify our problem. Much more attention is paid to it in \cite{RefMh} where the marks form radii of circles centered at the points of a point process given by a density with respect to the Poisson process. The random set corresponding to the union of circles in a compact window is investigated. Since an exact method is not available the authors use an approximation what means that estimation of the distribution of radii is done by methods for a Boolean model. Then an MCMC maximum likelihood method (\cite{RefMw}) is used for the estimation of parameters of the point process. A recent paper by \cite{RefDe} deals with the same model as \cite{RefMh}, their goal is to use the Takacs-Fiksel estimator instead of the computationally demanding maximum likelihood method. 

In our work we present a solution of the problem for another random set. We deal with the planar segment process (\cite{RefCh}, \cite{RefP}) having a density of exponential form with respect to a Poisson process. We consider first a model with reference directional distribution and in the end a model with reference length distribution of segments. The difference in comparison to marked point processes presented in the literature is that when the directional distribution is present we are not completely on Euclidean spaces. Our main tool is the derived relation between the reference and observed mark distribution for the segment process. The basic asymptotic properties of the Takacs-Fiksel method of estimation are known, see \cite{RefCd}. In the present paper we pay attention to the computing of semiparametric estimators using simulated data and quantifying the variability for small sample size.

First some background from spatial point processes having a density is presented. Then in Section 3 a general formula for the mark distributions is derived. In Section 4 we present a Gibbs type segment process with reference directional distribution and fixed lengths of segments. Interactions enter the model by means of intersections, cf. \cite{RefVb}, \cite{RefV}. An approximation from \cite{RefBn} is used to avoid the problem with unknown moments. The Takacs-Fiksel estimator for this model is developed in Section 5 and tested on data simulated by MCMC algorithm from \cite{RefG}. The estimator of the reference directional distribution is evaluated numerically. In Section 6 we address an old problem of the existence of a stationary process with given conditional intensity from the previous Section. Section 7 presents a model with reference length distribution of segments while the reference directional distribution is uniform. This model is an inhomogeneous Poisson process with a condition on segments to lie entirely in a circular window. The maximum likelihood method of estimation is available, see Section 8, we use the isotropy of the process to simplify the computation.

\section{Spatial point process given by a density}\label{PointProcess}
Consider a bounded Borel set $B\subset {\mathbb R}^d$ with Lebesgue measure
$|B|>0$ and a measurable space $({\mathbf N},{\mathcal N})$ of integer-valued
finite measures on $B.$ ${\mathcal N}$ is the smallest $\sigma $-algebra
which makes the mappings $\x\mapsto \x(A)$ measurable for all Borel sets
$A\subset B.$ A random element having a.s. values
in $({\mathbf N},{\mathcal N})$ is called a finite point process. Let a
Poisson point process $\eta $ on $B$ have finite intensity measure $\lambda $
with no atoms and distribution $P_\eta $ on $\mathcal N.$ We consider a
finite point process $X$ on $B$ given by a density $p$ w.r.t. $P_\eta ,$
i.\,e. with distribution $P_X $ 
\begin{equation}\label{dns}{\mathrm d}P_X
(\x)=p(\x)\,{\mathrm d}P_\eta (\x),\; \x\in \mathbf N,\end{equation} where
$p:{\mathbf N}\rightarrow {\mathbb R}_+$ is measurable satisfying
$$\int_{\mathbf N}p(\x)\,{\mathrm d}P_\eta (\x)=1.$$  As
described in \cite{RefB}, p.\,61, integer-valued finite measures can be
represented in this context by $n$-tuples of points corresponding to their
support ($n$ is variable). We have \begin{equation}\label{vzp}\p(X\in A)=\end{equation}$$=e^{-\lambda (B)}\sum_{n=0}^\infty\frac{1}{n!}\int_B\dots\int_B{\bf 1}[(x_1,\dots ,x_n)\in A]p(x_1,\dots ,x_n)\lambda (\dd x_1)\dots\lambda (\dd x_n),$$ $A\in\mathcal N.$ For an hereditary density $p$ the distribution of the process is alternatively determined by the conditional intensity $$\lambda^*(\x ,u)=\frac{p(\x\cup {u})}{p(\x)},\; \x\in\N ,\; u\notin\x.$$
An important tool is the Georgii-Nguyen-Zessin (GNZ) formula \begin{equation}\label{gnz}\E\left[\sum_{u\in X}q(u,X\setminus u)\right]=\int_B\E[\lambda^*(X,u)q(u,X)]\lambda (\dd u),\end{equation} valid for any measurable test function $q$ on $B \times {\mathbf N}.$

\section{Segment process with reference mark distribution}\label{SegmentProcess}
In the paper we study random segment processes in the plane $\R^2.$ A segment is a closed set which will be parametrized by its centre $y=(y_1,y_2)\in\R^2,$ length $r>0$ and direction $\varphi\in[0,\pi ).$ There is a bijection $\iota $ between the parametric space $\{u=(y,r,\varphi )\}$ and the system of segments as a subsystem of the space of closed sets in $\R^2.$ Throughout the paper we use exclusively the parametric representation of segments, omitting the bijection $\iota $ in some expressions, e.g. $z\in u$ means a point $z$ of the segment $u,$ etc.

A segment process can be considered as a marked point process with two marks corresponding to the length and direction of a segment. Let $B\subset \R^2$ be bounded measurable,
\begin{equation}\label{bck}Y=B\times S,\quad S=(0,b]\times [0,\pi ),\end{equation}
where $b>0$ is an upper bound for the segment length, $[0,\pi )$ is the manifold of axial directions. Further $({\mathbf N},{\mathcal N})$ is a measurable space of integer-valued
finite measures on $Y.$ Let the Poisson process
$\eta $ on $Y$ have intensity measure $\lambda ,$
\begin{equation}\label{inti}\lambda (\dd(y, r, \varphi ))=\frac{1}{b\pi} {\mathrm d}y\dd r\dd\varphi
,\end{equation} where on the right hand side we have a multiple of Lebesgue measure on $Y.$ Let the segment process $X $ have an hereditary density $p$ with respect to $\eta $ 
\begin{equation}\label{hust}p(\x)=c\tau^{n(\x )}\prod_{u\in\x}g(r,\varphi )\exp (bD(\x )+aE(\x )),\;\x\in{\mathbf N},\end{equation} where $c>0$ is a normalizing constant, $b\in\R,$ $a\leq 0,$ and $\tau >0$ are parameters, $g$ is the density of reference length-direction distribution, $n(\x)$ is number of segments in the configuration $\x .$ Further \begin{equation}\label{eneg}D(\x )=\sum_{u\in\x}d(u),\quad E(\x)=\sum^{\neq}_{\{u,v\}\subset\x}\Phi(u,v)\end{equation} (sum over pairs of different segments from $\x$) for some nonnegative functions $d,\;\Phi$ on $Y,\;Y^2,$ respectively, $\Phi$ is called the pair potencial, $E$ energy function in the theory of Gibbs processes. 
The conditional intensity $$ \lambda^*(\x ,u)=\tau g(r,\varphi )\exp\left(bd(u)+a\sum_{v\in\x\setminus\{u\}}\Phi(u,v)\right).$$ 

Let $\rho $ be the intensity function of the process $X,$ i.e. for $A\subset Y$ $$\E X(A)=\int_A\rho (u)\dd u.$$ 
Using $q(u,\x)={\bf 1}_{[u\in A]}$ in the GNZ formula we obtain \begin{equation}\label{podmi}\rho (u)=\E\lambda^*(X,u).\end{equation}  
For $A\subset B$ measurable and the point process $\nu $  of segment centres, denote $$\kappa (A)=\mathbb{E}\sum_{y\in\nu }1_A(y).$$ 
Let $D\subset S,$  $$C(A\times D)=\mathbb{E}\sum_{(y,\xi)\in X} 1_A(y)1_D(\xi ).$$ 
The measure $\omega_D,\;\omega_D(A)=C(A\times D),$ is absolutely continuous with respect to $\kappa $ and the Radon-Nikodym density $P^y(D),$ such that \begin{equation}\label{vetap}\omega_D(A)=\int_AP^y(D)\kappa (\dd y),\end{equation} is the Palm mark distribution of $X$ at $y.$ 
Let $f_{X}^{(y)}$ be the density of the Palm mark distribution of length and direction of a typical segment of the process $X $ at the location $y\in B.$ 
\begin{proposition} For each $y\in B$ we have \begin{equation}\label{hst} f_{X}^{(y)}(\xi )=\frac{\rho(y,\xi )}{\int_{S}\rho(y,\xi )\dd\xi }.\end{equation}
\end{proposition}{\bf Proof:} For Borel sets $A\subset B,\; D\subset S$ it holds using the Campbell theorem $$C(A\times D)=\int_{Y} 1_A(y)1_D(\xi)\rho(y,\xi )\dd\xi\dd y.$$ Specially for $D=S$ we have $$\kappa (A)=\int_A\int_{S}\rho(y,\xi )\dd\xi \dd y.$$ Then from (\ref{vetap}) we have
$$\int_AP^y(D)\int_{S}\rho(y,\xi)\dd\xi\dd y=\int_A\int_D\rho(y,\xi)\dd\xi \dd y$$ for all Borel sets $A\subset B,$ finally
$$P^y(D)=\frac{\int_D\rho(y,\xi)\dd\xi}{\int_{S}\rho(y,\xi)\dd\xi }$$
 and (\ref{hst}) follows.
\hfill $\Box $
\section{Segment process with reference directional distribution}\label{RefDirec}
In this Section we assume that the segment length $r>0$ is fixed and for a bounded Borel set $B\subset {\mathbb R}^2$ we deal with the 
\begin{equation}\label{bck1}Y=B\times [0,\pi ).\end{equation} A segment $u=(y,\varphi )\in Y$ has center $y$ and axial orientation $\varphi .$ 
Consider a measurable space $({\mathbf N},{\mathcal N})$ of integer-valued
finite measures on $Y,$ alternatively of the supports of these measures.
 
We deal with the unit Poisson segment process $\eta $ with the intensity measure $$\lambda (\dd y)=\dd y\frac{1}{\pi }\dd\varphi $$ on $Y.$ 
Let the segment process $X$ have a density $p$ with respect to $\eta ,$ we consider \begin{equation}\label{den}p(\x)=c\exp(a\, N(\x))\tau^{n(\x )}\prod_{u_i\in\x}g(\varphi_i),\;\x\in {\mathbf N}, \end{equation}
where $N(\x )$ is the total number of intersections between segments, $g$ a reference probability density on $[0,\pi ),$ $\varphi_i$ direction of $i-$th segment $u_i,$ 
$a \leq 0,\; \tau >0$ are parameters, $c$ a normalizing constant. That means in the general model (\ref{hust}) we have $b=0,$ $g$ depends only on $\varphi $ and the pair potential is \begin{equation}\label{pote}\Phi (u,v)={\bf 1}_{[u\uparrow v]},\end{equation} i.e. indicator of the event that segment $u$ hits segment $v.$ This model belongs to the more general class of facet processes (\cite{RefVe}). The conditional intensity is \begin{equation}\label{papan}\lambda^*(\x ,u)=\tau g(\varphi )\exp (aN_\x (u)),\end{equation} where $N_\x (u)$ is the number of segments of $\x\setminus\{u\}$ hit by the segment $u.$

{\bf Example:} Let the direction density $g$ be that of von Mises distribution on $[0,\pi )$ with parameters $\kappa\geq 0,\;\mu\in \R ,$ $$g(\phi )=h(\kappa )\exp (\kappa\cos (2(\phi -\mu ))),\;\phi\in [0,\pi ],\quad h(\kappa )=\frac{1}{\pi I_0(\kappa )},$$ $I_0(\kappa )$ is the modified Bessel function of the first kind and order 0. Then \begin{equation}
\label{dens}p(\x)=c(\theta )\exp(\langle\theta ,\,G(\x)\rangle),\;\x\in {\mathbf N},\end{equation}
where $$\theta =(a,\,\log (\tau h(\kappa )),\,\kappa),$$$$G(\x)=(N(\x ),\,n(\x),\,\sum_{{x_i}\in\x}\cos (2(\varphi_i-\mu ))),$$ segments $x_i=(u_i,\varphi_i)$ have centres $u_i\in B$ and directions $\varphi_i\in [0,\pi ).$ 
The normalizing constant $c(\theta )$ is defined as
\[ c(\theta)= \frac{1}{ \int \exp\{\langle\theta,\,G(\x)\rangle \} {\mathrm d}P_{\eta } (\x)} .\] 
Also $$\{\theta\in {\mathbb R}^3:\; \int \exp\{\langle\theta,\,G(\x)\rangle \} {\mathrm d}P_{\eta } (\x)<\infty\}$$ is the largest set of $\theta $ such that the density (\ref{dens}) is well defined. 

In the present paper we want to relax the unimodality assumption of the reference directional distribution and deal with a general density $g.$
Let $f_X^{(y)}$ be the density of the Palm mark distribution $P^y$ of direction of a typical segment of the process $X $ at the location $y\in B.$ For $u=(y,\varphi )$ we get from (\ref{papan}) that \begin{equation}\label{aprx}\rho (u)=\tau g(\varphi )\E_X\exp (aN_\x (u))\end{equation} is the intensity function of $X,$ cf. (\ref{podmi}). The expectation $\E_X$ with respect to the distribution of $X$ is not analytically tractable, therefore
\cite{RefBn} suggest an approximation \begin{equation}\label{apx}\E_X e^{aN_\x (u)}\approx\E_{\eta (\rho )} e^{aN_\x (u)},\end{equation} where $\eta (\rho )$ is a Poisson process with intensity function $\rho .$
\begin{proposition}For $u=(y,\varphi )\in Y$
we have $$\E_{\eta (\rho )} e^{aN_\x (u)}=\exp\left((e^{a}-1)\int_{M_u}\rho(x)\dd x\right),$$ where $M_u=\{x\in Y;\; x\uparrow u\}$.\end{proposition} {\bf Proof:} Using (\ref{vzp}) we have $$\E_{\eta (\rho )} e^{aN_\x (u)}=$$$$=\exp \left(-\int_Y\rho (x)\dd x \right )\sum_{n=0}^{\infty }\frac{1}{n!}\int_Y\dots \int_Y\exp (aN_{(x_1,\dots ,x_n)}(u))\rho (x_1)\dd x_1\dots \rho (x_n)\dd x_n$$$$=e^{-\int_Y\rho (x)\dd x}\sum_{n=0}^{\infty }\frac{1}{n!}\prod_{i=1}^n\int_Y\exp (a 1_{[u\uparrow x_i]})\rho (x_i)\dd x_i=\exp\left(\int_Y(e^{a1_{[u\uparrow x]}}-1)\rho (x)\dd x\right)=$$$$=\exp\left((e^{a}-1)\int_{M_u}\rho(x)\dd x\right).$$ \hfill $\Box $

From Proposition 1 we have $$\rho (y,\beta )=C_yf_X^{(y)}(\beta ),\;\beta\in [0,\pi ),$$ for normalizing constants $C_y,\; y\in B.$  Assuming that there exists a stationary segment process in $\R^2$ with given conditional intensity (this time the conditional intensity cannot be defined by means of densities, but from the energy function) corresponding to (\ref{papan}), we have that $C=C_y,\;f_X^{(y)}=f_X$ do not depend on $y$ in the extension of $B$ onto the whole $\R^2.$ Under this assumption, discussed in Section 6, we estimate ($\beta $ is the direction of segment $x$) $$\int_{M_u}f_X(\beta )\dd x\approx\int_{\R^2\times [0,\pi )}{\bf 1}_{[x\uparrow u]}f_X(\beta )\dd x=r^2\int_0^\pi |\sin (\beta-\varphi)|f_X(\beta )\dd\beta .$$
 Using (\ref{aprx}) we can then express the desired density $g$ approximately as
\begin{equation}\label{gfi}g(\varphi)\approx\frac{Cf_X(\varphi )}{\tau\exp ((e^a-1)Cr^2J(\varphi ))},\end{equation} where $$J(\varphi )=\int_0^\pi |\sin(\varphi -\beta)|f_X(\beta )\dd\beta .$$ 

\section{Semiparametric estimation using Takacs-Fiksel approach}\label{TakacsFiksel}
In this section we suggest a method of estimation of parameters $C,a,\tau $ and the density $g$ from the previous section using the Takacs-Fiksel method. From formula (\ref{gnz}) we obtain innovation equations
\begin{equation}\label{tri}\delta_\x (q)=\sum_{u\in X}q(u,X\setminus u)-\int_Y\lambda^*(X,u )q(u,X)\dd u=0\end{equation} and solve them for various test functions $q.$ We take $\lambda^*$ from (\ref{papan})
where we insert approximation (\ref{gfi}) for unknown $g.$ First the density $f_X$ is estimated using a kernel estimator for directional data \cite{RefM}. Then put
$$\beta (a,r,C,\varphi )=\exp ((e^a-1)r^2CJ(\varphi )),$$ and we estimate
$C,a$ from the system of Takacs-Fiksel equations:
$$\sum_{u \in \x} N_{\mathbf{x}}(u) =\frac{\pi |B|C}{J}\sum_{i=1}^J\frac{f_{X}(\varphi_i)N_\x (u_i) e^{a N_{\mathbf{x}}(u_i)}}{\beta (a,r,C,\varphi_i )},$$ 
$$n(\mathbf{x})=\frac{\pi |B|C}{J}\sum_{i=1}^J\frac{f_{X}(\varphi_i)e^{a N_{\mathbf{x}}(u_i)}}{\beta (a,r,C,\varphi_i )}.$$ Here in the innovations equations we take score functions $q(u,\x)=N_\x(u)$ and $q=1,$ respectively, the integrals in the second term of (\ref{tri}) are evaluated by Monte Carlo method using $J$ independent simulations of segments uniformly distributed in $Y.$ Then we plug the estimators of $C,a$ in a formula obtained by integrating (\ref{gfi}):
$$\tau =\frac{\pi C}{J}\sum_{i=1}^J\frac{f_{X}(\varphi_i)}{\beta (a,r,C,\varphi_i )}$$
and finally estimate $g$ from (\ref{gfi}). 

A numerical study is based on twice 100 simulated realizations of segment process with parameters $\kappa = 1, \mu=0, \tau= 1000, r= 0.12$ on $[0,1]^2\times [-\frac{\pi }{2},\frac{\pi }{2}]$. The two cases I, II investigated are $a = -0.5,\; a=-3,$ respectively. The results are in Table \ref{tbb1}, Fig. 1 and Fig. 2. 

In Table 1 we observe a small difference between the true and mean values for both estimates of $a$ and $\tau $. The coefficient of variation $$CV=\frac{\rm sd}{|{\rm mean}|}$$ is also comparable despites the fact that the model II involves more interactions (inhibition of intersections) than the model I. In Fig. 1 we can observe how the kernel estimator of the observed (Palm mark) directional distribution differs from the true reference directional distribution (von Mises). The results in Fig. 2 suggest that the estimate of the reference density is slightly better (smaller bias and variability) for the case I than for the case II. We conclude that the approximation (\ref{apx}) works well in the Takacs-Fiksel method here.
\begin{table}
\centering
\begin{tabular}{ccccc}
\hline
 & true & mean & sd & CV  \\
\hline
 $a$ & -0.5 & -0.496 & 0.071 & 0.14 \\ 
$\tau $ & 1000 & 1011 & 154.7 & 0.15 \\
\hline
& true & mean & sd & CV \\
\hline
 $a$ & -3 & -3.03 & 0.356 & 0.12 \\ 
$\tau $ & 1000 & 976 & 141.0 & 0.14 \\
\hline
 \end{tabular}
\caption{Empirical mean, standard deviation (sd) and coefficient of variation (CV) of Takacs-Fiksel estimates of scalar parameters in the model having density (\ref{den}) with reference directional distribution. It is based on 100 simulations, the two cases correspond to $a=-0.5,\; a=-3.$}
\label{tbb1}
\end{table}
\begin{figure}
\begin{center}
\includegraphics[width=5.6cm, height=5cm]{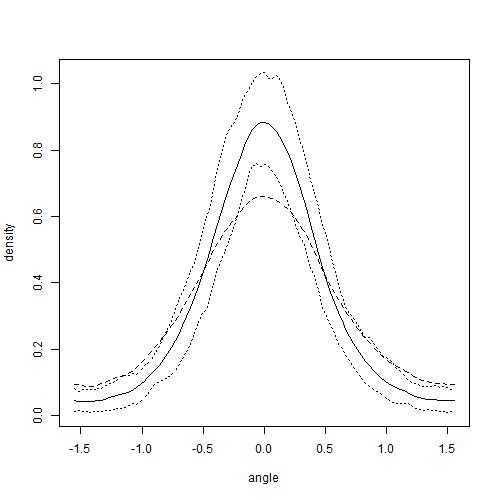}
\includegraphics[width=5.6cm, height=5cm]{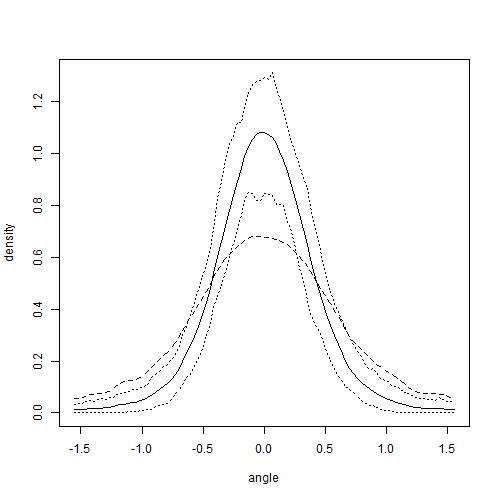}
\caption{Kernel estimation of the observed directional density based on 100 simulations of the segment process $X,\; a=-0.5$ (left), $a=-3$ (right). The average kernel estimator of the observed directional density (full line) compared to the true reference density (dashed line) of von Mises distribution with parameters $\mu =0,\;\kappa =1.$ The envelopes (dotted lines) correspond to empirical 90\% confidence interval for the kernel estimator, pointwise in 100 points on horizontal axis.}
\end{center}
\end{figure}
\begin{figure}
\begin{center}
\includegraphics[width=5.6cm, height=5cm]{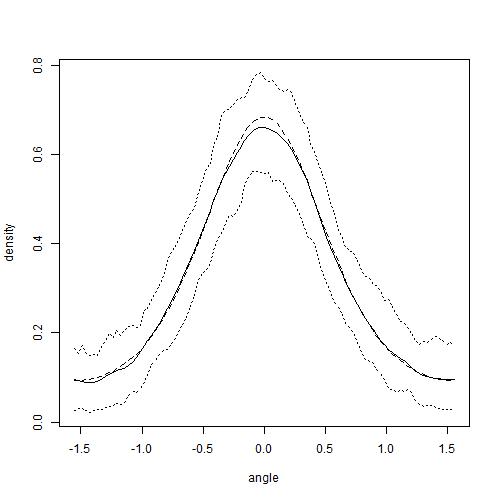}
\includegraphics[width=5.6cm, height=5cm]{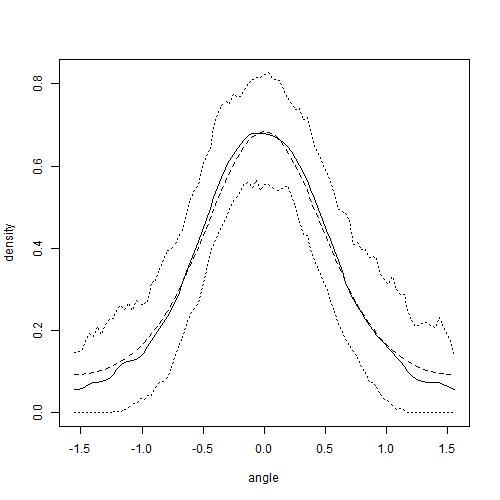}
\caption{Semiparametric estimation based on 100 simulations of the segment process $X,\; a=-0.5$ (left), $a=-3$ (right). The average estimator of the reference density (full line) compared to the true reference density (dashed line) of von Mises distribution with parameters $\nu =0,\;\kappa =1.$ The envelopes (dotted lines) correspond to empirical 90\% confidence interval for the estimated reference density, pointwise in 100 points on horizontal axis.}
\end{center}
\label{trih}
\end{figure}
\section{Existence of a stationary process with given conditional intensity}\label{ExistStac}
In order to be able to use the approximation (\ref{gfi}) correctly we need a sufficient condition for the existence of a stationary Gibbs segment process in $\R^2$ with prescribed conditional intensity. Various conditions are present in the literature starting with \cite{RefR}. We shall use a recent work of \cite{RefD} who deals with the concept of an hereditary energy function, invariant with respect to shifts, satisfying the finite range property. While his results are formulated for point processes in $\R^d,$ what we need here from \cite{RefD} is valid also for particle processes in the sense of \cite{RefSW}. This straightforward extension of a part of \cite{RefD}from Gibbs point process to Gibbs particle process is presented in \cite{RefNB}, specially also for segment processes in the plane. The energy function is hereditary if for each $\x\in{\mathbf N},\;y\in supp\,\x$ we have $$H(\x)<\infty\implies H(\x\setminus\{y\})<\infty ,$$ where $supp$ is the support of a measure.
\begin{definition} A function $f:{\mathbf N}\rightarrow \R$ is local (on $A$) if there exists $A\subset \R^d$ bounded such that for all $\x\in{\mathbf N}$ it holds $$f(\x)=f(\x_A),$$where $\x_A$ is the restriction of $\x$ onto $A.$\end{definition}
\begin{definition} An energy function $H$ has finite range $R>0$ if for all $A\subset\R^d$ bounded the local energy $$H_A(\x)=H(\x)-H(\x_{A^c})$$ is a local function on $A\oplus b(0,R).$\end{definition}
\begin{proposition}The energy function $E$ from (\ref{eneg}) with pair potential (\ref{pote}) has finite range property with $R=r.$\end{proposition}
{\bf Proof:} Let $A\subset\R^d$ be bounded, $\x\in{\mathbf N}.$ We have to show that $$E_A(\x)=E_A(\x_{A\oplus b(0,r)}),$$ i.e. $$E(\x)-E(\x_{A^c})=E(\x_{A\oplus b(0,r)})-E(\x_{(A\oplus b(0,r))\cap A^c}).$$ This is true since to both sides exactly intersections of such pairs of segments contribute, which have either both centres in $A$ or one centre in $A$ and the other in $(A\oplus b(0,r))\cap A^c.$ \hfill $\Box $ 
\begin{corollary} There exists a stationary segment process in $\R^2$ with conditional intensity (\ref{papan}).\end{corollary}{\bf Proof:} The energy function $E$ from (\ref{eneg}) with pair potential (\ref{pote}) is nonnegative, invariant with respect to shifts, hereditary and has finite range. According to \cite{RefD} and \cite{RefNB} the existence of a stationary segment process is guaranteed. \hfill $\Box $

The uniqueness issue is more complex, see \cite{RefD} for Gibbs point process in $\R^d.$ Here the uniqueness is not investigated, we deal with an existing stationary process in order to use the approximation (\ref{gfi}).

\section{The segment process with reference length distribution}\label{RefLength}
Consider a circle $B\subset {\mathbb R}^2$ centered in the origin with diameter $e_a>0.$ Let $L_o=[0,e_a]$ be the interval of segment lengths. Then
\begin{equation}\label{bck1}Y=B\times L_o\times [0,\pi )\end{equation} is the space of segments $u=(y,r,\varphi )\in Y$ which have centre $y,$ length $r$ and axial direction $\varphi .$  

We consider the Poisson segment process $\eta $ with the intensity measure $$\lambda (\dd u)=\dd y\frac{1}{e_a}\dd r\frac{1}{\pi }\dd\varphi $$ on $Y.$ 
Let the segment process $X$ have a density $p$ w.r.t. $\eta :$ \begin{equation}\label{den1}p(\x)=c1_{[\x\subset B]}\exp(b\, D(\x))\tau^{n(\x )}\prod_{u_i\in\x}f_1(r_i), \end{equation} $b\in\R,\,\tau >0,\; c$ is the normalizing constant, $r_i$ is the length of $i-$th segment $u_i,\;  f_1$ is a reference probability density on $L_o$ and
$$D(\x)=\sum_{u\in\x }d(u),\quad d(u)=\max_{z\in u}\frac{||z||}{e_a}.$$ That means in the general model (\ref{hust}) we have $a=0,$ $g=f_1$ depends only on the length variable $r,$ marginal reference directional distribution is uniform. Moreover a factor is added in (\ref{den1}) which forces the segments to lie completely in the window. 
 
For $u\subset B$ the quantity $d(u)\in (0,\frac{1}{2}]$ is the normalized distance of the most distant point of $u$ from the centre of $B.$ For positive, negative values of parameter $b$ more, less distant segments from the origin prevail, respectively, cf. Fig. 3. 

\begin{figure}
\begin{center}
\includegraphics[width=12cm]{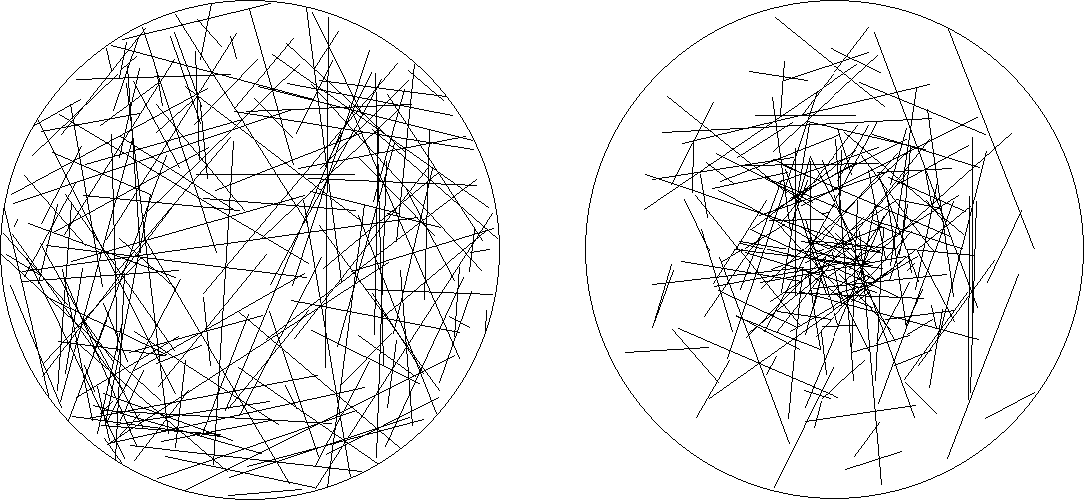}
\caption{Simulated realizations of the process with $e_a=1$ and the density (\ref{den1}) having parameters $b=10,\;\tau=3$ (left) and $b=-10,\;\tau=4000$ (right), $f_1$ corresponds to beta distribution with parameters $\alpha =2,\;\beta =4.$}\end{center}
\label{siml}
\end{figure}

The corresponding conditional intensity is
$$\lambda^*(\x ,u)=1_{[\x\cup\{u\}\subset B]}\tau f_1(r)\exp(bd(u)).$$ For the intensity function it holds 
$$\begin{array}{ll} \rho (u)=\E\lambda^*(X,u)=\tau f_1(r)\exp(bd(u)), & u\subset B,\\
\rho (u)=0, & u\cap (\R^2\setminus B)\neq\emptyset . \end{array}$$
There are no interactions among the segments in the model (\ref{den1})   with statistics $D(\x ),$ $X$ is in fact an inhomogeneous Poisson process with unknown reference density $f_1$ and a condition $X\subset B.$  The process $X$ is isotropic since both the reference Poisson process and the density $p$ are invariant with respect to rotations around the origin. We denote $Ay$ rotation of $y\in\R^2,$ where $A$ is an orthonormal matrix $2\times 2,$ analougously $Au$ is rotation of the segment $u$ and finally we use the same symbol $A\varphi $ for the corresponding direction of the segment $u$ after rotation $Au.$  Then also the intensity function is invariant with respect to rotations, i.e. $$\rho (u)=\rho (Au)$$ for all $u\in Y$ and all rotations $A.$ 

 Let $f_X^{(y)}$ be the bivariate (Palm) density of the distribution of length and direction of the segment centered at $y,$ from Proposition 1 we have for $y\in B$ a normalizing constant $C_y>0$ such that \begin{equation}\label{delk}f_1(r)=\frac{C_yf_X^{(y)}(r,\varphi )}{e^{bd(u)}}.\end{equation} From the isotropy of $X$ it holds \begin{equation}\label{rott}f_X^{(y)}(r,\varphi )=f_X^{(Ay)}(r,A\varphi )\end{equation} and therefore $$C_y=C_{Ay}$$ for all $y\in B$ and all rotations $A.$

\section{Semiparametric estimation using the maximum likelihood method}\label{PseudoLikel}
The aim is to estimate parameters $b,\tau $ and density $f_1$ from simulated data. We are using again a semiparametric approach so that $f_1$ is not parametrized. Because of inhomogeneity of $X$ the solution  has to be discretized, but we make use of isotropy. In the parametric part, since the process is Poisson we use maximum likelihood method for parameter estimation. For an observed realization $\x $ the likelihood is defined as $${\cal L}=\prod_{u\in\x}\lambda^*(\x \setminus u,u)\exp\left(-\int_{v\subset B}\lambda^*(\x ,v)\dd v\right).$$ The logarithmic likelihood $$\log {\cal L}=\log (\tau )n(\x)+bD(\x)+\sum_{u\in\x}\log f_1(r)-\tau\int_{u\subset B}f_1(r)e^{bd(u)}\dd u$$ has to be maximized with respect to $\tau, b.$ We have $$\frac{\partial\log {\cal L}}{\partial\tau}=\frac{n(\x )}{\tau }-\int_{u\subset B}f_1(r)e^{bd(u)}\dd u,$$$$\frac{\partial\log {\cal L}}{\partial b}=D(\x )-\tau\int_{u\subset B}d(u)f_1(r)e^{bd(u)}\dd u.$$ Using (\ref{delk}) we obtain equations \begin{equation}\label{ml}n(\x )=\tau\int_{u\subset B}C_yf_X^{(y)}(r,\varphi )\dd u,\end{equation}$$D(\x )=\tau\int_{u\subset B}C_yd(u)f_X^{(y)}(r,\varphi )\dd u.$$ 
In the estimation procedure we proceed in several steps:

\noindent (i) Consider discrete levels of $||y||$ like $0<w_1<\dots <w_k<e_a/2$ and put $y_j=(0,w_j)^T\in B.$
For $\Delta_f=\frac{e_a}{2k}$ let $$w_j=(j-1)\Delta_f+\frac{\Delta_f}{2},\; j=1,\dots ,k,$$
$$Y_j=B_j\times L_o\times [0,\pi ),\quad B_j=\{y\in B,\,(j-1)\Delta_f<||y||\leq j\Delta_f\},\;j=1,\dots ,k.$$ 
The kernel estimator of the bivariate densities $f_X^{(y_j)},\;j=1,\dots ,k,$ is evaluated from the observed data in each class, i.e. from the sample of segments $u_j^{(i)}=(y_j^{(i)},r_j^{(i)},\varphi_j^{(i)}),\; i=1,\dots ,m_j$ centered in $B_j,\; j=1,\dots ,k.$ Because of (\ref{rott}) these segments are first rotated by such $A_j^{(i)}$ that $$A_j^{(i)}y_j^{(i)}=qy_j$$ for some $q>0.$ Then we apply kernel estimation to each sample $$A_j^{(i)}u_j^{(i)}, i=1,\dots , m_j, $$ to estimate the length-direction density $f_X^{(y_j)},\; j=1,\dots ,k.$ Since the length component has values on a compact $L_o,$ we use here a system of beta kernels \cite{RefC}. The angular component is a circular variable, cf. \cite{RefM}, jointly we use a bivariate product kernel.

\noindent (ii) 
Since the length of the longest possible segment centered in $j-$th class is $$l_j=2\sqrt{\left(\frac{e_a}{2}\right)^2-(\Delta_f(j-1))^2},$$ we use the step $\Delta_c^{(j)}=l_j/m$ for numerical integration of (\ref{delk}) with $y=y_j$ and a fixed $\varphi $ using Simpson rule, to express unknown constants $C_{y_j}=C_j,\; j=1,\dots ,k,$ as functions of $b.$

\noindent (iii) Integrals in the equations (\ref{ml}) are evaluated by Monte Carlo simulation of $M$ segments in $Y$ uniformly randomly. Let $n$ of them lie completely in $B,$ denoted $\bar{u}_j^{(i)}=(\bar{y}_j^{(i)},\bar{r}_j^{(i)},\bar{\varphi}_j^{(i)}),\; i=1,\dots ,n_j$ where $n_j$ of them are centered in $B_j,\;\sum_{j=1}^k n_j=n.$ Then \begin{equation}\label{btu}\frac{D(\x)}{n(\x)}=\frac{\sum_{j=1}^k C_j \sum_{i=1}^{n_j}d(\bar{u}_{j}^{(i)})f_X^{(\bar{y}_{j}^{(i)})}(\bar{r}_j^{(i)},\bar{\varphi}_j^{(i)})}{\sum_{j=1}^k C_j \sum_{i=1}^{n_j}f_X^{(\bar{y}_{j}^{(i)})}(\bar{r}_j^{(i)},\bar{\varphi}_j^{(i)})}\end{equation} is an equation with a single variable $b$ which is solved numerically.

\noindent (iv) Having estimated $b$ we obtain $C_{j}$ from step (ii) and then $\tau $ from any of the equations
$$\tau=\frac{4M\,n(\x )}{\pi^2e_a^3\sum_{j=1}^k C_j\sum_{i=1}^{n_j}f_X^{(\bar{y}_{j}^{(i)})}(\bar{r}_j^{(i)},\bar{\varphi}_j^{(i)})},$$  $$\tau=\frac{4M\,D(\x )}{\pi^2e_a^3\sum_{j=1}^k C_j\sum_{i=1}^{n_j}d(\bar{u}_j^{(i)})f_X^{(\bar{y}_{j}^{(i)})}(\bar{r}_j^{(i)},\bar{\varphi}_j^{(i)})}$$ which are Monte Carlo analogues of the equations in (\ref{ml}). 

\noindent (v) Finally the estimator of the reference length density $f_1$ is obtained by plugging these estimators in (\ref{delk}) and renormalizing. We choose $\varphi =0,$ from the levels of $y$ those close to the boundary are mostly variable and therefore not used. We can average results of several levels of $y.$

\begin{table}
\centering
\begin{tabular}{ccccc}
\hline
 & true & mean & sd & CV \\
\hline
 $b$ & 3 & 3.051 & 0.481 & 0.16 \\ 
$\tau $ & 900 & 948.9 & 186.2 & 0.19 \\
\hline
\end{tabular}
\caption{Empirical mean, standard deviation (sd) and coefficient of variation (CV) of the maximum likelihood estimates of scalar parameters in the model (\ref{den1}) with reference length density. Based on a sample of 60 simulated realizations of the segment process with parameters $b=3, \tau= 900,\alpha =2,\beta =4, e_a=1.$}
\end{table}

A numerical study is based on a sample of 60 simulated independent realizations of the segment process with parameters $b=3, \tau= 900,\alpha =2,\beta =4, e_a=1.$ In step (i) of the estimation procedure $k=6$ classes were considered. In Fig. 4 there are kernel estimates of the observed length density in all six classes separately. We can observe the difference between the estimated observed length distribution and the true reference beta distribution, which increases along classes towards the boundary of the window. The results of estimation of parameters $b, \tau $ (step (iii), (iv)) are in Table 2. We observe still a reasonably small difference between the true and the mean value of estimates of $b$ and $\tau ,$ but it is larger than in the Gibbs model with directional distribution in Section 5, also CV is slightly larger. This can be justified by the fact that the present model is inhomogeneous and the estimation procedure is more complex. The semiparametric estimator of the reference length density $f_1$ (step (v)) is in 
Fig. 5, here we do not consider two outer classes where the inhomogeneity is the highest (since there are little admissible directions of longer segments close to the border). We observe a small bias of the estimator of the reference density, but the variability is larger than in the first model in Section 5 from the same reasons as we argued to Table 2.

\begin{figure}
\begin{center}
\includegraphics[width=12cm, height=8cm]{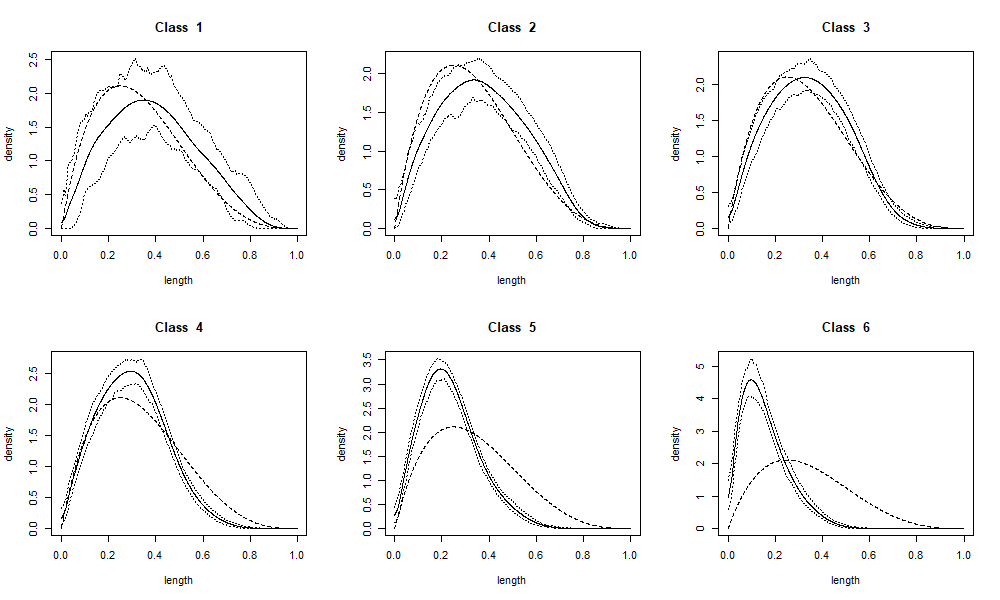}
\caption{Kernel estimation of the observed length density based on 60 simulations of the inhomogeneous segment process $X,\; b=3, \tau=900.$ In each of six classes the average kernel estimator of the observed length density (full line) is compared to the true reference density (dashed line) of beta distribution with parameters $\alpha =2,\;\beta =4.$ The envelopes (dotted lines) correspond to empirical 90\% confidence interval for the kernel estimator, pointwise in 100 points on horizontal axis.}
\end{center}
\end{figure}
\begin{figure}
\begin{center}
\includegraphics[width=5.6cm, height=5cm]{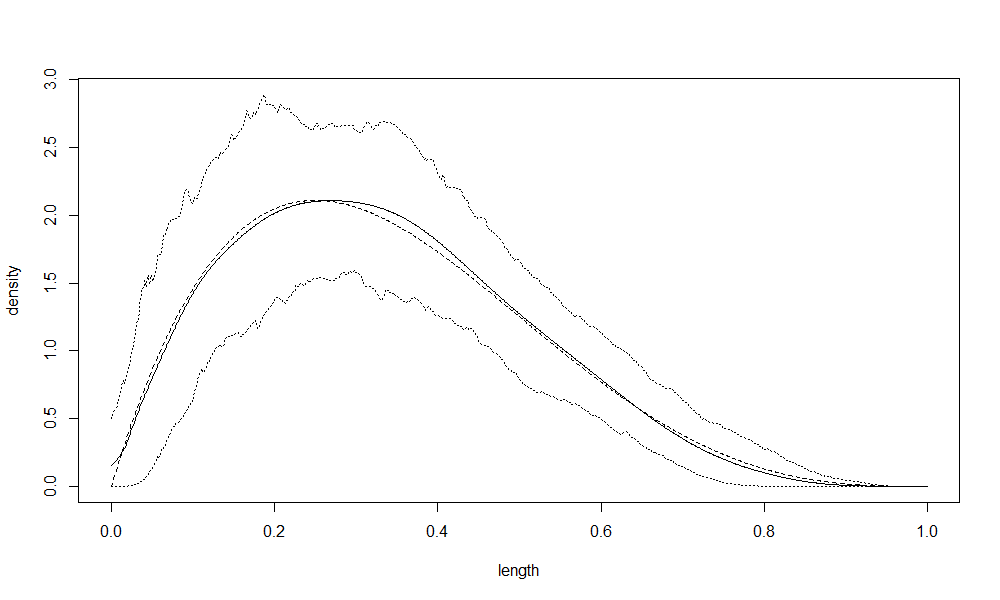}
\caption{Semiparametric estimation of reference length density based on 60 simulations, respectively, of the segment process $X,\; b=3,\tau=900.$ The average (from first four classes) estimator of the reference density (full line) compared to the true reference density (dashed line) of beta distribution with parameters $\alpha =2,\;\beta =4.$ The envelopes (dotted lines) correspond to empirical 90\% confidence interval for the estimated reference density, pointwise in 100 points on horizontal axis.}
\label{simla}
\end{center}
\end{figure}

\medskip

\subsection*{Acknowledgements} The research was supported by the Czech Science Foundation, project 16-03708S and by Charles University, grant SVV-2016-260334.


\bigskip

\begin{thebibliography}{}
\bibitem{RefBa} Baddeley A, Turner R (2000) Practical maximum pseudolikelihood for spatial point processes. Austr NZ J Statist 42:283--322.

\bibitem{RefB} Baddeley A (2007) Spatial point processes and their applications. In: Stochastic geometry. Ed by Weil W, Lecture Notes in Math, vol 1892, Springer, Berlin, 1--75.

\bibitem{RefBn} Baddeley A, Nair G (2012) Fast approximation of the intensity of Gibbs point processes. Electron J Statist
6:1155--1169. 

\bibitem{RefBe} Bene\v{s} V, Ve\v{c}e\v{r}a J, Eltzner B, Wollnik C, Rehfeldt F, Kr\'{a}lov\'{a} V, Huckemann S (2017) Estimation of parameters in a planar segment process with a biological application. Image Anal Stereol 36:25--33. 

\bibitem{RefC} Chen SX (1999) Beta kernel estimators for density functions. Comput Statist Data Anal 31:131--145.

\bibitem{RefCh} Chiu B, Stoyan D, Kendall WS, Mecke J. (2013) Stochastic Geometry and Its Applications. 3rd Ed, Wiley, New York.

\bibitem{RefCd} Coeurjolly JF, Dereudre D, Drouilhet R, Lavancier F (2012) Takacs-Fiksel Method for Stationary Marked Gibbs Point Processes. Scand J Statist 39:416--443.

\bibitem{RefD} Dereudre D (2017) Introduction to the theory of Gibbs point processes. Preprint arXiv:1701.08105 [math.PR], to appear in Lecture Notes, Springer. 

\bibitem{RefDe} Dereudre D, Lavancier F, Helisova K (2014) Estimation of the intensity parameter
of the germ-grain Quermass-interaction model when the number of germs is not
observed. Scand J Statist 41:809--829.

\bibitem{RefG} Geyer CJ, M{\o}ller J (1994) Simulation procedures and likelihood
inference for spatial point processes. Scand J Statist, 21:359--373.

\bibitem{RefM} Mardia KV, Jupp PE (1999) Directional Statistics. Wiley, New York.

\bibitem{RefMh} M{\o}ller J, Helisova K (2010) Likelihood inference for unions of interacting discs. Scand J Statist 37:365--381.

\bibitem{RefMw} M{\o}ller J, Waagepetersen R (2004) Statistical Inference and Simulation for Spatial Point Processes. CRC Press, Singapore.

\bibitem{RefNB} Novotn\'{a} D, Bene\v{s} V (2017) Central limit theorem for functionals of Gibbs
particle processes. Preprint arXiv:1707.06872 [math.PR], submitted to Kybernetika. 

\bibitem{RefP} Pawlas Z (2014) Self-crossing points of a line segment process. Methodol Comp Appl Probab 16:295--309. 

\bibitem{RefR} Ruelle D (1970) Superstable interactions in classical statistical mechanics. Comm Math Phys 18:127--159.

\bibitem{RefSW} Schneider R, Weil W (2008) Stochastic and Integral Geometry. Springer, Berlin.

\bibitem{RefV} Ve\v{c}e\v{r}a J (2016) Central limit theorem for Gibbsian U-statistics of facet processes. Appl Math 61, 4:423--441.

\bibitem{RefVb} Ve\v{c}e\v{r}a J, Bene\v{s} V (2016) Interaction processes for unions of facets, the asymptotic behaviour with increasing intensity. Methodol Comp Appl Probab  18, 4:1217--1239. 

\bibitem{RefVe} Ve\v{c}e\v{r}a J, Bene\v{s} V (2017) Approaches to asymptotics for $U$-statistics of Gibbs facet processes. Statist Probab Lett 122:51--57.
\end{thebibliography}
\end{document}